\documentstyle[12pt]{article}
\font\BBm=msbm10   scaled\magstep1

\newcommand{\beq}{\begin{equation}}
\newcommand{\eeq}{\end{equation}}
\newcommand{\be}{\begin{equation}}
\newcommand{\ee}{\end{equation}}

\def\ba{\begin{array}}
\def\ea{\end{array}}
\def\l{\label}
\def\ra{\rightarrow}

\def\Z{\BBm Z}         
\def\Z+{{\BBm Z}^0}    

\def\da{a^{\dag}}

\def\db{b^{\dag}}
\def\A{${\cal A}(q)$ }
\def\H0{${\cal H}_0$ }
\def\Hg{${\cal H}_{\gamma}$ }

\def\ZPC{{\it Z. Phys. C.}}
\def\JMP{{\it J. Math. Phys.}}

\begin{document}

\begin{flushright} {\large \bf FTUV/96-50,\,\,IFIC/96-58} \\

\vspace*{0.1cm}

{\begin{tabular}{c} q-alg/9610002 \\
September 1996 \end{tabular}}
\end{flushright}
\vspace*{.8cm}

\begin{center}
{\Large \bf IRREDUCIBLE REPRESENTATIONS} \\   \vspace*{2mm}
{\Large \bf OF DEFORMED OSCILLATOR ALGEBRA} \\ \vspace*{2mm}
{\Large \bf AND $q$-SPECIAL FUNCTIONS}
\end{center}

\vspace{.4cm}

\begin{center}
{\bf E. V. Damaskinsky} \\ \vspace*{1mm}
{\it Defense Constructing Engineering Institute, \\
Zacharievskaya st 22, 
191194, St.Petersburg, Russia}  \\ \vspace*{3mm}
{\bf P.P. Kulish} \footnote{Permanent address: St.Petersburg Branch of
the Steklov Mathematical Institute, Russian Academy of Sciences,
Fontanka 27, St.Petersburg 191011, Russia} \\ \vspace*{1mm}
{\it Departamento de F{\'{\i}}sica The{\'o}rica, Univ. de Valencia,\\
46100, Burjassot (Valencia),  Spain.}
\end{center}

\vspace*{.4cm}

\begin{center} {\bf Abstract} \end{center}
{\parindent=.8in
\narrower
Different generators of a deformed oscillator algebra
give rise to one-parameter families of $q$-exponential
functions and  $q$-Hermite polynomials related by 
generating functions. Connections of
the Stieltjes and Hamburger classical moment problems with
the corresponding resolution of unity for the $q$-coherent states 
and with 'coordinate' operators - Jacobi matrices, are also pointed out.
\par}

\newpage
{\large \bf 1.}
Suggestions to change the canonical commutation
relations for improving some properties of quantum
field theory have already appeared in the works of the founders
of quantum mechanics (see \cite{a1,FdP} and Refs therein).

In the early fifties, E. P. Wigner \cite{a1} posed the
following question: ``what kind of functions, $f(H)$,
appearing in the right-hand side of the commutator
$ \lbrack p,x]=i\hbar f(H)$
are compatible with the given expression of the Hamiltonian, $H$,
and the standard equations of motions''?  He found that, in the
case of the harmonic oscillator Hamiltonian 
the usual
expression $f(H)=\bf{1}$ was not unique. These investigations
were continued in \cite{a2} where, for a special case, now
called ``$q$ - root of unity'', it was found that, together
with the fermionic (two-dimensional) and bosonic (infinite
dimensional) cases, there exist cases with
dimensions equal to $m$ which are
related with the parastatistics not connected to the Green's
Ansatz.

Twenty years later, the generalization of the
Veneziano amplitude, by substitution of the $q$-$\Gamma$-function
instead of the standard $\Gamma$-function \cite{coon}, 
gave rise in the operator
formalism to the $q$-oscillator commutation relation \cite{AC},
\be \l{ac}  aa^{\dagger}-qa^{\dagger}a=1, \quad 0<q<1. \ee

A rejuvenation of the problem of oscillator deformations
in the end of the 80's was associated with the growing interest
in quantum groups. Such  popularity appeared after the
works \cite{B,M}, in which a deformation based on the relation
$AA^{\dagger}-q^{1/2}A^{\dagger}A=q^{-N/2}$, 
was considered in connection with Schwinger's realization of
the quantum algebra $su_q(2)$ \cite{KR} (for a $q$-boson
description of the $su_q(1,1)$, see \cite{KD}).
Further interest in the  $q$-oscillator problem was
stimulated by research in the multimode
case \cite{PW}, supersymmetries \cite{CKL}, and relations to
the $q$-analysis \cite{FV} (for more details and
Refs, see \cite{DK}).

Different generators of the deformed oscillator algebra give
rise to one-parameter families of $q$-exponential functions
\cite{DK,Ex,McA,At}, $q$-Hermite polynomials and 
other $q$-special functions \cite{DK2,Kl,Koor,NG}.
Consideration of the resolution of unity (completeness of the
system of $q$-coherent states) for the $q$-Bargmann - Fock
realization of irreducible representations of deformed
oscillator algebra, and the spectral properties \cite{Kl} of
the 'coordinate' operator (which is represented as
a Jacobi matrix), pointed out deep connections with the
classical  Stieltjes and Hamburger moment problems \cite{AcBe}.

Recently, the $q$-oscillator was applied to the study of the
phonon spectrum in ${}^4\!$He \cite{He}, a specific case of
the one-dimensional Schr\"odinger equation \cite{And},
different quantum mechanical models \cite{Wess}, and the trapped atom
problem. 

\vspace*{.5cm}
{\large \bf 2.} The deformed oscillator  algebra, \A, is generated
by three elements $a,\, \da , \, N$ with defining relations
\be\l{1} a \da-q \da a=1\, ,\quad 
[N,\, a]=-a, \quad [N,\, \da ]=\da \, . \ee
The generator $N$ is considered as an independent element, and
we restrict ourselves to the case of positive real
$q\in (0,\infty)$. The algebra \A has a central element \cite{PK},
\be\l{3} \zeta = q^{-N}\bigl( [N;q]-\da a \bigr)
;\,\quad [N;q]:= (1-q^N)/(1-q) \ee
(for a more general three-generator algebra \A with 
$[a,\,\da]=F(N)$, see \cite{Curt,KQ}).

In the original papers, the irreducible representation of \A with
the vacuum state $|0\rangle$ $(a|0\rangle = 0)$ was
considered. The oscillator-type representation space \H0, in
the basis of eigenvectors of the operator $N$, is
\be \l{4}  {\cal H}_0 = \{\,|n\rangle; \quad n=0,1,2,...;
\quad a|0\rangle=0, \quad |n\rangle
= ([n;q]!)^{-1/2}(a^{\dagger})^n|0\rangle\, \}. \ee

Due to the existence of a non-trivial central element,
$\zeta$, in addition to \H0, the algebra \A has a set of
inequivalent irreducible representations $(0<q<1)$ in the
spaces \Hg $(\gamma \geq {\gamma}_c = (1-q)^{-1})$
parameterised by the value of the central element
$\zeta=-\gamma$ \cite{PK}, with the spectrum of $N$, the set 
of all integers {\BBm Z}. The matrix $\da$ in the 
number operator basis is
\be\l{5} (\da)_{n\, k}=c_n \delta_{n\, k+1} \; , \quad \da
|n-1\rangle =c_n |n\rangle \;,\quad (c_n)^2=
\gamma q^n + [n;q] \; .\ee 
These irreducible representations are connected with different 
symplectic leaves of Poisson brackets in $R^3$ , which correspond 
to the quasiclassical limit of the $q$-oscillator commutation 
relation \cite{KD}.

Considering \A as an associative algebra, any invertible
transformation of the generators is admissible; in particular,
there are some natural sets of the generators: 
\be\l{6} A A^{\dagger}- q^{1/2}A^{\dagger}A =q^{-N/2} \, ,
\quad [N, A]=-A \, , \quad [N, A^{\dagger}]= A^{\dagger}, \ee
related to the quantum algebra $sl_q(2)$ via the Schwinger
realization \cite{B,M}, and the following set related to the
$sl_q(2)$ algebra by a contraction procedure with fixed $q$ \cite{PK},
\be\l{7} [\alpha, \alpha^{\dagger}]=q^{-N}, \quad [N,\alpha]=
-\alpha, \quad [N,\alpha^{\dagger}]=\alpha^{\dagger} . \ee
The equivalence of these generators is given by the equalities
$a=q^{N/2}\alpha = q^{N/4}A$ \cite{KD,DK}, with an obvious 
one-parameter generalization, namely, 
\be\l{8} a(\lambda)=q^{-{\frac12}\lambda N}a\, , \quad
\da(\lambda)=\da q^{-{\frac12}\lambda N}\, . \ee
This leads to the commutation relations (still one
degree of freedom)
\be\l{9} a(\lambda) \da(\lambda) - q^{1-\lambda} 
\da(\lambda) a(\lambda) = q^{-\lambda N}\, .\ee
Sometimes, these generators and relation (\ref{9}) are called a two-parameter
deformed oscillator \cite{2p} : $p$ $\leftrightarrow$ $q^{1-\lambda}$ and $r$ 
$\leftrightarrow$ $q^{-\lambda}\, ,$ $ a \da - p \da a = r^{N}\,.$ 
However, they define the same algebra \A in the case of general $q=p/r$.

One more formal parameter $\nu \in {\BBm R}$ can be added by
a shift $N \ra N+\nu$. The corresponding set of \A
generators is denoted by $W_{p,r}^{\nu} (q)$ \cite{BDY}.
As a consequence of (\ref{9}), namely, 
\be\l{11} a(\lambda) \bigl( \da(\lambda)\bigr) ^m =
\bigl( p \da(\lambda) \bigr) ^m a(\lambda)  + \bigl( p \da
(\lambda) \bigr) ^{m-1} r^{N} [m;{\frac rp}] \, , \ee
the normalized  basis vectors of \H0 in terms of $ \da(\lambda)$ are given by 
$$ |n\rangle = \bigl( [n;q,\lambda ]!\bigr)^{-1/2}
\bigl( \da(\lambda) \bigr) ^n |0\rangle$$  
with the factorials defined as
\be \l{13} [n;q,\lambda ]!=\prod\limits_{k=1}^{n}
[k;q,\lambda ]\, ,\quad [m;q,\lambda ] = q^{\lambda (1-m)}[m;q] \,. 
\ee

\vspace*{.5cm}

{\large \bf 3.} In the theory of Lie groups and quantum
mechanics, special functions appear as particular matrix
elements (overlap coefficients) of appropriate operators
in corresponding representations (realizations): examples
are exponential functions, as coherent states in the 
Bargmann-Fock representation of \H0, of the usual boson oscillator
$[b,\,\db]=1$,
\be\l{3.1} {\rm exp}({\overline{\!w}}z)=\langle w|z\rangle ,
\quad |z\rangle = {\rm e}^{z\db}|0\rangle ,
\quad b|z\rangle =z|z\rangle ,\ee
and Hermite polynomials, as eigenvectors of the operator $N$,
in the coordinate representation,
$$H_n(x) \sim \langle n|x\rangle,
     \quad (b+\db)|x\rangle=2x|x\rangle .$$ 
The simple action of the annihilation and creation operators in the 
coherent state representation leads to the generating function of 
the Hermite polynomials 
\be\l{gf} \omega (z;\,x) = \langle \bar{z}|x \rangle = 
 {\rm exp}(2xz -{\frac 12}z^2) \, .
\ee

The coherent states of the annihilation operator $a$ 
of the $q$-oscillator in \H0 were introduced in \cite{AC}:
\be\l{3.2} a |z\rangle =z|z\rangle , \qquad
|z\rangle ={\rm e}_q(z\da)|0\rangle \, , \ee
\be\l{3.3} {\rm e}_q(x)=\sum\limits_{n=0}^{\infty}
\frac{x^n}{[n;\,q]\,!} \, . \ee

In the $q$-Bargman-Fock space, related to the $q$-coherent
states, the creation operator $\da$ is the operator of multiplication
by $\bar{z}$,
\be\l{3.4} \ba{c} |f\rangle \ra f(z)=
   \langle {z}|f \rangle , \\[.2cm]
   \langle {z}|\da |f\rangle =
   (a\ |z\rangle )^{\dag}\,|f\rangle=
   \bar{z}\,\langle {z}|f \rangle =\bar{z}\,f(z),  \ea\ee
the annihilation operator, $a$, is a $q$-difference operator
${\cal D}_q$, and the $q$-coherent state is given by the 
$q$-exponent (\ref{3.3}),
\be\l{3.5}  a\,f(z) = {\cal D}_q\,f(z)
=\frac{f(z)-f(qz)}{z(1-q)}, \qquad
\langle z|\zeta \rangle  = {\rm e}_q(\bar{z}\zeta). \ee
The $q$-exponent above (see (\ref{3.3})) is well known in 
$q$-analysis \cite{Ex}. The scalar product in the 
$q$-Bargman-Fock realization of \H0 is given by \cite{AC}
\be\l{3.6} \langle \phi|f \rangle  =
     \frac{1}{2\pi}\int \overline{\phi(z)}\,
      f(z)\, {\rm d}\mu (z), \ee
where the measure is defined by the resolution of unity,
\be\l{3.7} \frac{1}{2\pi} \int_0^{1/(1-q)}\!\!\!
     \int_0^{2\pi}|z\rangle \langle z| \;\bigl(
     {\rm e}_q(q|z|^2) \bigr)^{-1}\,{\rm d}
     \phi\;{\rm d}_q|z|^2  =\sum_{n=0}^\infty
     |n \rangle \langle n| =I: \ee
this completeness relation was proved in \cite{AC} using
the product representation of the $q$-exponent (\ref{3.3}), namely,
\be \l{3.8}  {\rm e}_q(x) =
    \Bigg( \prod_{k=0}^\infty (1-(1-q)q^kx) \Bigg)^{-1} = 
\frac{1}{((1-q)x;q)_\infty}, \ee
and the Jackson $q$-integral, $\int_0^{b} f(x)\,{\rm d}_qx=(1-q) 
\sum_{m=0}^\infty q^{m}b\,f(q^{m}b) $ \cite{Ex}.

Following the same pattern, other choices for the generators of the
$q$-oscillator algebra \A give rise to different $q$-exponential
functions \cite{KD,DK},
\be \l{3.10} \alpha |z\rangle_\alpha = z|z\rangle_\alpha , 
\quad |z\rangle_\alpha  =
         {\rm e}_{1/q}(z\alpha ^{\dag}) |0\rangle,  \ee
\be \l{3.11} A |z\rangle_A = z|z\rangle_A , \quad |z\rangle_A
   = {\rm E}_{q}(zA^{\dag}) |0\rangle,	    \ee
where the symmetric $q$-exponent is 
\be \l{3.12} {\rm E}_{q}(x)=\sum_{m=0}^{\infty}
      \frac{x^m}{[m]_q\,!}, \qquad   [m]_q=
      \frac{q^m-q^{-m}}{q-q^{-1}}. \ee

The one-parameter $q$-exponential function 
${\rm exp}(z;q,\lambda)$ 
is connected with the annihilation operator $a(\lambda)$ (\ref{8}), (\ref{9}) 
\be \l{3.13} a(\lambda) |z;\lambda\rangle = 
    z|z;\lambda\rangle  , \quad |z;\lambda\rangle = 
   {\rm exp}(za(\lambda)^{\dag};q,\lambda) |0\rangle ;\ee
\be \l{3.14} {\rm exp}(z;q,\lambda) =\sum_{m=0}^{\infty}
    q^{\lambda \,{n\,(n-1)}/2}\frac{x^n}{[n,q]\,!}.  \ee
The properties of these $q$-exponents  
$({\rm exp}(z;q,\lambda))$ are quite different 
\cite{McA,At}; for example, for $0<q<1$ and $\lambda < 
0$, the $q$-exponent  
${\rm exp}(z;q,\lambda)$ (\ref{3.14}) has zero radius of
convergence. It would be interesting to relate different 
$q$-exponential functions and their properties with particular 
physical systems. 

The corresponding resolution of unity in the 
$(q,\lambda)$-Bargman-Fock realization of \H0, where the annihilation 
operator (\ref{8}) acts as a difference operator ${\cal D}_q^{(\lambda)} $ 
\cite{McA}, 
\be\l{3.15} \frac{1}{2\pi} \int_0^{\infty}\!\!\!
     \int_0^{2\pi}|z;\lambda\rangle \langle z;\lambda| \;
     {\rm d}\phi \;{\rm d}_q\sigma(|z|^2) =
     \sum_{n=0}^\infty |n \rangle \langle n| =I,  \ee
results in the classical moment problem (MP) \cite{AcBe} 
for the measure ${\rm d}_q\sigma(|z|^2)$,
\be\l{3.16} \int_0^{\infty}\,x^n\,{\rm d}_q\sigma(|z|^2) 
       = s_n(q;\lambda),  \quad s_n(q;\lambda) =
       [n;q,\lambda]\,!\, .   \ee

Depending on the behaviour of the moments $s_n(q;\lambda)$ as 
$n\ra \infty$, the MP can be determinate (a unique solution, 
if any: this is the case of the $q$-oscillator 
(\ref{1})), or indeterminate (many solutions: these cases
are realized for the $q$-oscillators (\ref{6}) or (\ref{7})). The 
completeness (the system is overcomplete) was proved for 
$s_n(q;\lambda)=[n;q]\,!$ \cite{AC}, 
$s_n(q;\lambda)=[n]_q\,!$ \cite{Br}, and  
$s_n(q;\lambda)=[n;q^{-1}]\,!$ \cite{PPK}. Complete subsystems 
of $q$-coherent states (\ref{3.2}) (or (\ref{3.13}) for 
$\lambda=0$) are discussed in \cite{Per}.

The classical MP refers also to $q$-Hermite polynomials: the 
latter are nothing but polynomials of the first kind 
\cite{AcBe} for a Jacobi matrix ${\cal J}$ which is 
constructed as a ``generalized coordinate'' from the 
$q$-oscillator creation and annihilation operators 
\cite{DK2},
\be\l{3.17}  {\cal J}(\lambda)= 
a(\lambda) + \da(\lambda),\quad {\cal J}(\lambda)\, 
     |x\rangle_{\lambda} = 2x \,|x\rangle_{\lambda}, \ee
\be\l{3.18}  |x\rangle_{\lambda}
     =\sum_{n=0}^{\infty}H_n (x;q,\lambda)|n\rangle . \ee
Due to (\ref{3.17}), these $q$-Hermite polynomials satisfy the
following three-term recurrence relation:
\be\l{3.19} c_n(\lambda)H_{n-1}(x;q,\lambda) +
      c_{n+1}(\lambda)H_{n+1}(x;q,\lambda)=  
      x\,H_n(x;q,\lambda). \ee 
The corresponding generating function can be introduced as in the 
oscillator case (\ref{gf}), $\omega(z,x;\lambda) = 
\langle \bar{z} ;\lambda|x\rangle_{\lambda}, $
however its form will depend on the chosen 
generators of \A  \cite{DK2} 
$$\langle \bar{z} ;\lambda|( a(\lambda) + \da(\lambda) )|x\rangle_{\lambda} = 
({\cal D}_q^{(\lambda)} + z) \, \omega(z,x;\lambda) = 
2x \, \omega(z,x;\lambda) \,.$$ 
This difference equation for $\omega(z,x;\lambda) $ will include two points 
for $\lambda = 0,\,1$, so its solution will be given by the 'standard' 
$q$-exponent (\ref{3.3}) and three points: $z,\, q^{-\lambda}z,\, 
q^{1-\lambda}z $ for the general $\lambda$. 
The measure entering into the $q$-Hermite polynomials 
$H_n (x;q,\lambda)$ orthogonality 
relations is connected with the solution of the Hamburger MP: 
this measure is known explicitly for some cases
(see e.g.\cite{Kl}). This connection of the MP with Jacobi matrices 
gives rise to a generalized deformation of the oscillator 
identifying the matrix $c_k \delta_{n+1,k}\, , c_k > 0$ with an annihilation 
operator $a$. Then one gets the Wigner commutation relation 
$[a,\, \da] = F(N)$ with $F(n) = c_{n+1}^2 - c_n^2$ and its central 
element $\zeta = \bigl(c^2(N) -\da a \bigr) + const $
(see also \cite{FdP,Curt,KQ}). The $q$-special functions related
to the other irreducible representations \Hg of \A are discussed 
in \cite{PPK}. In particular, for the generators (2) the normalized 
$q$-coherent states exist in \Hg for the creation operator $\da $
and $z > \gamma_c = (1-q)^{-1}$. 

\vspace*{.5cm}

{\bf Acknowledgements.} The authors express their gratitude to the 
organizers of the workshop IWCQIS-96 (JINR, Dubna) and especially 
to G. Pogosyan for 
hospitality. The authors thank J.A. de Azcarraga and A. Vartanyan for useful 
discussions. This research was supported by the RFFI 
grants N 95-01-00569-a (EVD) and N 96-01-00851 (PPK). The second 
author thanks the Generalitat Valenciana for financial support.


\begin{thebibliography}{99}
\small{

\bibitem{a1}  Wigner E.P., {\it Phys.Rev.}, {\bf 77}, 711-715 (1950)

\vspace{-8pt}

\bibitem{FdP}  Rampacher H., Stumpf H., Wagner F.,
        {\it Fortschr. Phys.}, {\bf 13}, 385 (1965)

\vspace{-8pt}

\bibitem{a2}  O'Raifeartaigh L., Ryan C.,
    {\it Proc. Roy. Irish. Ac.}, ser.A., {\bf 62}, 93-115 (1963)\\
    Gruber B., O'Raifeartaigh L.,
    {\it Proc. Roy. Irish. Ac.}, ser.A., {\bf 63}, 69-73 (1964)

\vspace{-8pt}

\bibitem{coon} Coon D.D., Baker M.,
        {\it Phys. Rev. D.}, {\bf D2}, 2349 (1970); \\
        Coon D.D., Yu S., Baker M.M.,
        {\it Phys. Rev. D.}, {\bf D5}, 1429 (1972);

\vspace{-8pt}

\bibitem{AC} {\rm Arik M., Coon D.D.},
             {\it J.  Math. Phys.}, {\bf 17}, no.4, 524-527 (1976).

\vspace{-8pt}

\bibitem{B}  Biedenharne L.C.,
              {\it J. Phys. A.}, {\bf 22}, no.18, L~873-878 (1989)

\vspace{-8pt}

\bibitem{M} Macfarlane A.J.,
     {\it J. Phys. A.}, {\bf 22}, no.21, 4581-4586 (1989)

\vspace{-8pt}

\bibitem{KR}  Kulish P.P., Reshetikhin N.Yu.,
           {\it Zap. Nauch. Semin. LOMI}, {\bf 101}, 101 (1981);

\vspace{-8pt}

\bibitem{KD} Kulish P.P., Damaskinsky E.V.,
           {\it J. Phys. A.}, {\bf 23}, no 9, L~415-419 (1990)

\vspace{-8pt}

\bibitem{PW} Pusz W., Woronowicz S.L.,
     {\it Repts Math. Phys.}, {\bf 27}, no.2, 231-257 (1989)

\vspace{-8pt}

\bibitem{CKL} Chaichian M., Kulish P.P.,
     {\sl Phys.Lett.B.,} {\bf 234},  no.1/2, 72-80 (1990);
Pr-t. CERN-TH 5969/90, \\
     Chaichian M., Kulish P.P., Lukierski J.,
     {\it Phys. Lett. B.}, {\bf 262}, no.1, 43-48 (1991).

\vspace{-8pt}

\bibitem{FV}  Floreanini R., Vinet L.,
    {\it Lett. Math. Phys.}, {\bf 22}, no.1, 45-54 (1991);
    {\it Phys. Lett. A.}, {\bf 170}, no.1, 21-28 (1992);
    {\it Ann. Phys.}, (N.Y.), {\bf 221}, no.1, 53-70 (1993).

\vspace{-8pt}

\bibitem{DK} Damaskinsky E.V., Kulish P.P.,
      {\it Zap.Nauch.Semin. LOMI}, {\bf 189}, 37-74 (1991) (in Russian)
       \\ English transl: {\it J. Soviet. Math.} {\bf 62}, 2963 (1992);\\
   Damaskinsky E.V., Kulish P.P.,
   {\it Quantum groups, deformed oscillators and their interrelations,}
   Proc. Workshop {\it ``Finite dimensional integrable systems''},
   JINR, Dubna, July, 1994

\vspace{-8pt}

\bibitem{Ex} Exton H., {\it q-hypergeometric functions and  applications},
    Chichester. Ellis Horwood.  1983;\\
    Gasper G., Rahman M., {\it Basic Hypergeometric Series},
    (CUP, Cambridge, 1990)

\vspace{-8pt}

\bibitem{McA} McAnally D.,
      {\it J. Math. Pphys.}, {\bf 36}, no.1, 546-573; 574-595 (1995)

\vspace{-8pt}

\bibitem{At}  Atakishiev N.M.,
     {\it J. Phys.A :Math.Gen.} {\bf 29}, no.10, L~223-227 (1996)

\vspace{-8pt}

\bibitem{DK2}  Damaskinsky E.V., Kulish P.P.,
    {\it Zap. Nauch. Sem. POMI,\/} {\bf 199}, 81-90 (1992) (in Russian);

\vspace{-8pt}

\bibitem{Kl} Burban I.M., Klymik A.U.,
     {\it Lett. Math. Phys.}, {\bf 29}, no.1, 13-18 (1993);\\
     Chung W.-S., Klimyk A.U.,
     {\it J. Math. Phys.}, {\bf 37}, no.2, 917-932 (1996)

\vspace{-8pt}

\bibitem{Koor} Koornwinder T.H., {\it Special functions and $q$-commuting 
     variables}, q-alg/9608008, 1996.

\vspace{-8pt}

\bibitem{NG} Nelson~C.A., Gartley~M.G., {\it On the two $q$-analogue 
      logarithmic functions}, q-alg/9608015, 1996.

\vspace{-8pt}

\bibitem{AcBe} Akhiezer N.I., {\it The classical moment problem},
	  Olyver and Boyd, Edinburgh, 1985 \\
           Berezanskii Ju.M.,
	  {\it Eigenfunction expansion of self-adjoint operators,}  AMS, 1968.

\vspace{-8pt}

\bibitem{He}  R-Mointeiro M., Rodrigues L.M.C.S., Wulck S.,
        {\it Phys. Rev. Lett.}, {\bf 76}, no.7, 1098-1101 (1996)

\vspace{-8pt}

\bibitem{And} Andrianov~A.A., Cannata~F., Dedonder~J.-P., Ioffe~M.V.,
     {\sl Phys. Lett. A.}, {\bf 217}, 7-14 (1996) 

\vspace{-8pt}

\bibitem{Wess} Lorek~A., Wess~J., \ZPC, {\bf 67}, no.4, 671-679 (1995)

\vspace{-8pt}

%

\bibitem{PK}  Kulish P.P., {\it Teor. Math. Phys.}
              {\bf 85}, no 1, 158-161(1991)

\vspace{-8pt}

\bibitem{Curt}  Curtright T.L.,
      in: {\it Quantum Groups}, 
(T.L.~Curtright, D.~Fairlie and C.K.~Zachos (Eds.)),
           Proceedings of the Argonne Workshop
          (World Scientific, Singapore, 1990).

\vspace{-8pt}

\bibitem{KQ} Katriel J., Quesne C., \JMP , {\bf 37}, no~4, 1650-1661 (1996)

\vspace{-8pt}

\bibitem{2p}  Chakrabarti R., Jagannatham R.,
       {\it J. Phys. A.}, {\bf 24}, no.13, L~711-718 (1991)


\vspace{-8pt}

\bibitem{BDY} Borsov V.V., Damaskinsky E.V., Yegorov S.B.,
     Stockholm pr-t TRITA-MAT-1995-MA-20; {q-alg/9509022}.

\vspace{-8pt}

\bibitem{Br} Bracken A.J., McAnally D.S., Zhang R.B., Gould~M.D.,
         {\it J. Phys. A.}, {\bf 24}, no~7, 1379-1392 (1991).

\vspace{-8pt}

\bibitem{PPK} Kulish P.P., {\it Irreducible representations of deformed
        oscillator and coherent states}, Preprint KTH-96/21, 11pp.,
        Stockholm, 1996.

\vspace{-8pt}

\bibitem{Per} A.M. Perelomov,
           {\it Helv. Phys. Acta}, {\bf 68}, 554--576 (1996).



}
\end{thebibliography}
\end{document}